\documentclass{llncs}
\usepackage[caption=false]{subfig}

\RequirePackage{lmodern}
\RequirePackage{graphicx}
\RequirePackage{epsfig}
\RequirePackage{amssymb, amsmath}
\RequirePackage{subfig}
\RequirePackage{stmaryrd}
\usepackage{physics,tabularx, environ, stmaryrd}
\usepackage{float}
\usepackage{pgfplots}
\usepackage{pgfplotstable}
\usepackage{tcolorbox}
\usepackage[lined,boxed,commentsnumbered, ruled,vlined,linesnumbered]{algorithm2e}
\pgfplotsset{compat=1.18}
\usepackage{tikz}
\usetikzlibrary{quantikz2}
\usetikzlibrary{plotmarks}
\usetikzlibrary{shapes.geometric}
\usetikzlibrary{arrows.meta}

\begin{document}

\title{Quantum model for CVRPTW}

\def\shorttitle{Titre court}

\author{Imran Meghazi\inst{1,2}, Éric Bourreau\inst{2}}

\institute{
    La Poste, France \\
    \and
    LIRMM, France \\
    \email{\{imran.meghazi, eric.bourreau\}@lirmm.fr}
}

\maketitle


\begin{abstract}{}
    This paper proposes a quantum algorithm for the capacitated vehicle routing problem with time windows (CVRPTW) based on Grover Search framework.
    This problem is often faced by Postal services in the context of package delivery or other time-sensitive operations.
    We provide an implementation on gate based quantum computer of a model inspired by classical route first, cluster second technique. The quantum paradigm allows to overcome suboptimality inherent
    property of this decomposition. In the current NISQ (Noisy Intermediate-Scale Quantum) era, the most important limitation is the number of available qubits which makes time windows and capacity constraints hard to tackle.
    We introduce a qubit-efficient split-inspired modeling which adds only a linear number of decision qubits to standard quantum formulations for Traveling Salesman Problem (TSP).
\end{abstract}


\section*{Introduction}

\par Operations research aims to tackle a wide range of combinatorial problems, among which the Vehicle Routing Problem (VRP) is a prominent example. The VRP generalizes the well known Travelling Salesman Problem (TSP).
The TSP seeks the shortest cycle that visits every node in a graph exactly once. In a logistics context, it can be interpreted as finding the most efficient route for a single vehicle to deliver one item to each customer.
The VRP extends this formulation by considering multiple vehicles operating simultaneously, each subject to capacity constraints, in order to serve all customers while minimizing the total travel cost.
This problem arises daily in logistics companies and organizations, such as the French postal service \textit{La Poste}.
In practice, delivery services face additional constraints. In this paper, we focus on two of them : a capacity constraint, which reflects the limited space of each vehicle, and a time window constraint that allows customers to specify a preferred delivery period.
Solving the VRP with such constraints, referred to as the Capacitated Vehicle Routing Problem with Time Windows (CVRPTW), is computationally difficult and is usually adressed using heuristics \cite{koskosidisOptimizationBasedHeuristicVehicle1992} and bio-inspired algorithms \cite{PDFApplyingANT}.
\par Quantum Computing is a promising paradigm that has already produced algorithms outperforming their classical counterpart. In some cases, such as Shor's algorithm \cite{shorPolynomialTimeAlgorithmsPrime1997} it offers exponential speedup,
in others, like Grover's algorithm \cite{groverFastQuantumMechanical1996}, a quadratic improvement. Grover's algorithm provides an oracle-based general framework for addressing combinatorial problems.
While quantum oracle-based models for the Traveling Salesman Problem and QUBO formulations for the CVRPTW have been proposed, to the best of our knowledge, no quantum oracle-based formulation has yet been introduced for the CVRPTW.
In this work, we propose such a model using Grover Adaptive Search \cite{durrQuantumAlgorithmFinding1999} on gate-based quantum computers.

\section{The Capacitated VRP with Time Windows}

\subsection{Problem definition}
A \textit{Capacitated Vehicle Routing Problem with Time Windows} can be defined as follows: a fleet of vehicles must
deliver packages to a set of $n$ customers $V = \{v_1, v_2, \ldots, v_n\}$, where all packages are initially stored at a central
depot $d$. Each customer $v_i$ requests a specific quantity of goods denoted by $q_i$. In addition to capacity constraints,
customers specify a preferred delivery time window during the day. In this work, we consider the case where each customer is assigned
a single time window, meaning that customer $v_i$ must be served within the interval $[a_i, b_i]$.

\subsection{Methods}
In classical operations research, the Vehicle Routing Problem is typically addressed using two main approaches.
The most well-known is the \textit{cluster-first, route-second} heuristic, which involves grouping customers into clusters
and then solving a separate routing problem, such as a TSP, within each cluster. For instance, customers may be
clustered based on their geographical proximity. It is then assumed that solving a TSP for each cluster yields a VRP solution
that is reasonably close to optimal. However, this assumption breaks down in the presence of time window constraints,
which can invalidate geographically efficient clusterings. Several heuristics based on this strategy can be found in the literature,
such as in \cite{koskosidisOptimizationBasedHeuristicVehicle1992}.

An alternative strategy, albeit less intuitive,  is the \textit{route-first, cluster-second} approach,
originally introduced by Beasley \cite{beasleyRouteFirstCluster1983}. In this method, a so-called \textit{giant tour}
is constructed that visits all customers, and this tour is subsequently divided into feasible \textit{sub-tours} that satisfy
the problem constraints. We specifically focus on the \textit{split} procedure presented in \cite{prinsSimpleEffectiveEvolutionary2004a},
where the partitioning of the giant tour is typically performed according to vehicle capacity limits.

\subsection{Split Procedure}

In this section, we detail the \textit{split} procedure, as our quantum approach is heavily inspired. 

The \textit{split} algorithm takes as input an auxiliary graph $H = (X, A, W)$ constructed from the solution obtained in the
first routing phase, namely the \textit{giant tour}. The set $X$ contains the nodes, ordered from node 0 to $n$ according to this tour.
Let $P$ be this tour. The set of arcs $A$ contains an arc $(i,j)$, with $i < j$, if a trip visiting every customers between $P_{i+1}$ and
$P_{j}$ doesn't violate any constraints. The corresponding weight $w_{ij}$ represents the total cost of that sub-tour, computed as
$D_{d, i+1} + \sum_{k=i+1}^{j-1} D_{P_k, P_{k+1}} + D_{j, d}$ with $D$ the distance matrix.
Finding a minimum cost path from node $0$ to $n$ in $H$ yields an optimal splitting of the giant tour, where every arc corresponds
to a feasible sub-tour that respects the original problem's constraints.

As an illustrative example, let us consider a typical instance of the Capacitated Vehicle Routing Problem.
\label{example}

\begin{minipage}[t]{0.5\textwidth}
    \centering
    \textbf{Distance matrix}
    \vspace{0.5em}

    \begin{tabular}{|c|*{7}{>{$}c<{$}|}}
        \hline
          & d  & 1  & 2  & 3  & 4  & 5  & 6  \\ \hline
        d & 0  & 23 & 30 & 23 & 14 & 20 & 26 \\ \hline
        1 & 23 & 0  & 17 & 27 & 27 & 38 & 36 \\ \hline
        2 & 30 & 17 & 0  & 21 & 40 & 46 & 37 \\ \hline
        3 & 23 & 27 & 21 & 0  & 35 & 40 & 12 \\ \hline
        4 & 14 & 27 & 40 & 35 & 0  & 16 & 31 \\ \hline
        5 & 20 & 38 & 46 & 40 & 16 & 0  & 33 \\ \hline
        6 & 26 & 36 & 32 & 12 & 31 & 33 & 0  \\ \hline
    \end{tabular}
\end{minipage}%
\begin{minipage}[t]{0.5\textwidth}
    \centering
    \textbf{Nodes}
    \vspace{0.5em}

    \begin{tabular}{|c|c|}
        \hline
        $v_i$ & \textbf{$q_i$} \\
        \hline
        $d$   & 0              \\
        $1$   & 2              \\
        $2$   & 3              \\
        $3$   & 1              \\
        $4$   & 3              \\
        $5$   & 2              \\
        $6$   & 3              \\
        \hline
    \end{tabular}

    $C^{max} = 5$
\end{minipage}

\vspace{0.5cm}
The only constraints is that, for every subtour, the total must not exceed the vehicle's maximum capacity $C^{max}$. Stated differently, the cumulative load at $k^{th}$ node, denoted $c_{P_k}$,
must not exceed $C^{max}$.
\begin{align}
    \sum_{k=i+1}^j c_{P_k} \leq C^{max}, ~\forall(i, j) \in X, i < j \label{eq:capconst}
\end{align}

\begin{figure}[H]
    \centering
    \subfloat[]{
        \centering
        \begin{tikzpicture}[scale=0.4]
    \node[draw=red, fill=red!20, rectangle, minimum size=1pt, font=\tiny] (D) at (0,0) {\textbf{d}};
    
    \foreach \x/\y/\num/\demande in {
        1.5/1.8/1/2, 2.8/1.2/2/3, 2.2/-0.8/3/1,
        -1.2/0.8/4/4, -1.8/-0.8/5/2, 1.2/-1.8/6/3
    } {
        \node[draw=blue, fill=blue!20, circle, minimum size=1pt, font=\tiny] (C\num) at (\x,\y) {\num};
    }
    
    \draw[red, -{Stealth[scale=0.9]}] 
        (D) -- (C1) -- (C2) -- (C3) -- (C6) -- (C5) -- (C4) -- (D) -- cycle;
\end{tikzpicture}
        \label{fig:giantour}
    }
    \subfloat[]{
        \centering
        \begin{tikzpicture}[scale=0.5, every node/.style={font=\tiny},
    arc/.style={->, font=\scriptsize, inner sep=1pt}
]
    \foreach \i/\label in {0/d, 1/1, 2/2, 3/3, 4/6, 5/5, 6/4} {
        \node[draw=blue!80!black, fill=blue!20, circle, minimum size=3pt] 
            (N\i) at (\i*1.5, 0) {\label};
    }

    \draw[->] (N0) -- node[above] {46} (N1);  
    \draw[->] (N1) -- node[above] {60} (N2);  
    \draw[->] (N2) -- node[above] {46} (N3);  
    \draw[->] (N3) -- node[above] {52} (N4);  
    \draw[->] (N4) -- node[above] {40} (N5);  
    \draw[->] (N5) -- node[above] {28} (N6);  

    \draw[red, arc, bend right=40] (N0) to node[below] {70} (N2);  

    \draw[arc, bend left=45] (N1) to node[above] {74} (N3);  
    \draw[red, arc, bend right=55] (N2) to node[below] {61} (N4); 
    \draw[arc, bend left=45] (N3) to node[above] {79} (N5);  
    \draw[red, arc, bend right=55] (N4) to node[below] {50} (N6); 
\end{tikzpicture}
        \label{fig:aux}
    }
    \subfloat[]{
        \centering
        \begin{tikzpicture}[scale=0.5, every node/.style={font=\tiny}]
    \node[draw=red, fill=red!20, rectangle, minimum size=3pt] (d) at (0,0) {\textbf{d}};
    

    \node[draw=blue, fill=blue!20, circle, minimum size=3pt] (C1) at (1.5,1.8) {1};

    \node[draw=blue, fill=blue!20, circle, minimum size=3pt] (C2) at (2.8,1.2) {2};

    \node[draw=blue, fill=blue!20, circle, minimum size=3pt] (C3) at (2.2,-0.8) {3};

    \node[draw=blue, fill=blue!20, circle, minimum size=3pt] (C4) at (-1.2,0.8) {4};

    \node[draw=blue, fill=blue!20, circle, minimum size=3pt] (C5) at (-1.8,-0.8) {5};

    \node[draw=blue, fill=blue!20, circle, minimum size=3pt] (C6) at (1.2,-1.8) {6};

    \draw[red, -{Stealth[scale=1.2]}] 
        (D) -- (C1) -- (C2)-- (D) -- (C3) -- (C6) -- (D) -- (C5) -- (C4) -- (D)-- cycle;
    
\end{tikzpicture}
        \label{fig:sol}
    }
    \caption{The steps to compute a solution. \ref{fig:giantour} the giant tour. \ref{fig:aux} the auxiliary graph $H$. \ref{fig:sol} the solution.}
\end{figure}
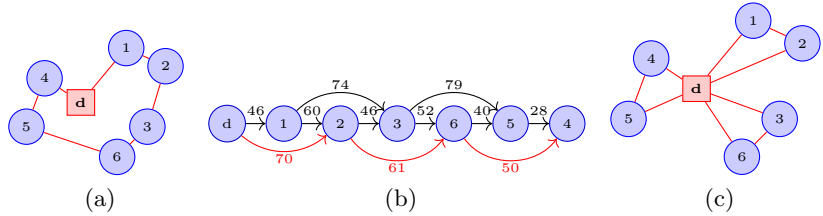

We first compute a giant tour, let $P$ be the ordering of this tour, $P_i$ is the $i^{th}$ node of the tour. Figure \ref{fig:giantour} show this giant tour and
figure \ref{fig:aux} shows the auxiliary graph.

\vspace{-0.5cm}

\par Finding an optimal split for this giant tour can be effectively done by computing the shortest path over figure \ref{fig:aux}.

\section{Quantum Model}

In classical computing, a split-based formulation is inherently suboptimal, as the optimal VRP solution is not necessarily
the optimal split of the optimal \textit{giant tour}. Consequently, it requires evaluating every possible \textit{giant tour} and
all its splits, an exponential task.
As a result, it is typically employed only as a highly efficient heuristic. In contrast, quantum superposition enables
the parallel evaluation of all these possibilities, allowing the identification of the optimal \textit{giant tour} for every possible
splits.

\subsection{Split Inspired model}
\label{sec:qsplit}

\paragraph{Objective function}
Let $G = (V, E)$ be the graph instance. Let $D$ denotes the distance matrix of this graph.
Following the approach introduced in the previous section, we define separately the decision variables related to the ordering of customers
and those related to the clustering (or splitting) of the route. Let $P$ represent the \textit{giant tour}, and let $P_i$ denote its $i^{\text{th}}$ node.

Then we define the following variables :
\[
    y_i =
    \begin{cases}
        1 & \text{if } P_i \text{ is the end of a sub-tour } \\
        0 & \text{otherwise}
    \end{cases}
\]

The variables $y_i$ encode the splitting decisions. Specifically, if $y_i = 1$, it indicates that the vehicle must return to the facility immediately
after visiting node $P_i$. Consequently, we enforce $y_0=1$ and $y_n = 1$ to ensure that the tour starts and ends at the depot.

Thus the value to minimize is defined as follows :
\begin{align}
    \sum_{i = 1}^{n-1} \left[D_{P_i, P_{i+1}} \left(1-P_i\right) + \left(D_{P_i, f} + D_{f, P_{i+1}}\right)y_i\right] + D_{f, P_1}+ D_{P_n, f}
    \label{eq:obj}
\end{align}

As discussed in the previous section, the clustering process (or splitting, in this context) must be performed in accordance with the problem's constraints.

\paragraph{Uniqueness constraint}
Each customer must be visited exactly once, meaning that no customer may appear more than once in the tour.
\begin{align}
    P_i \neq P_j, ~ \forall (i, j) \in V \label{eq:quniconst}
\end{align}

\paragraph{Capacity constraint}
This constraint ensures that the total load of each sub-tour does not exceed the vehicle’s maximum capacity.
Given a fixed customer ordering, it can be formalized as shown in Equation~\ref{eq:capconst}.
In practical terms, this means that at each node $i$, the cumulative load $c_i$ must be less than or equal to the vehicle's capacity limit $C^{max}$.

\begin{align}
    c_i \leq C^{max}, ~ \forall i \in V \label{eq:qcapconst}
\end{align}

In the quantum computing context, all operations must be reversible. As a result, it is necessary to store the load at each preceding node in order
to compute the current load. The load at node $i$ is thus calculated as the sum of the loads of the previous nodes, conditioned by the $y_i$ variables,
plus the load associated with node $i$ itself. As previously defined, $y_{i-1} = 1$ indicates that node $i-1$ marks the end of a sub-tour,
implying that the next node and --- thus next subtour--- is served by a new vehicle.

\begin{align}
    c_{i} = q_{P_i} + \left(c_{i-1} \times \left(1 - y_{i-1}\right) \right) \label{eq:qcapdef}
\end{align}

\paragraph{Time window constraint}
This constraint behaves similarly to the capacity constraint. The visiting time at node $i$, denoted by $t_i$, must not exceed the latest acceptable arrival time,
denoted by $b_i$.

\begin{align}
    t_{i} \leq b_{P_i}, ~ \forall i \in V \label{eq:qtimeconst}
\end{align}

As with the capacity constraint, it is necessary to store the visiting times of all previous nodes in order to compute the arrival time at node $i$.
In this case, when $y_i = 1$, the current time is reset, as each sub-tour is assumed to be performed by a separate vehicle.
Therefore, the computation of the time window constraint follows a similar logic to that of the capacity constraint. The key difference lies in the use of a maximum:
the arrival time at node $i$ is defined as the maximum between the earliest acceptable time $a_i$ and the cumulative travel time.

\begin{align}
    t_{i} = \max\left(a_{P_i}, T_{f, P_i} \times y_{i-1} + \left(\left(t_{i-1} + T_{P_{i-1}, P_i}\right)\times \left(1 - y_{i-1}\right) \right) \right) \label{eq:qtimedef}
\end{align}

This gives us the following model :

\begin{subequations}
    \begin{align}
        \min_{P, y}         & \sum_{i = 1}^{n-1} \left[D_{P_i, P_{i+1}} \left(1-y_i\right) + \left(D_{P_i, depot} + D_{depot, P_{i+1}}\right)y_i\right] + D_{depot, P_1}+ D_{P_n, depot} \nonumber \\*
        \textrm{s.t.} \quad & P_i \neq P_j, \quad \forall (i, j)                                                                                                                                   \\
                            & c_{i} \leq C^{max}, \quad \forall i                                                                                                                                  \\
                            & t_{i} \leq b_{P_i}, \quad  \forall i                                                                                                                                 \\
                            & c_{i} = q_{P_i} + \left(c_{i-1} \times \left(1 - y_i\right) \right)                                                                                                  \\
                            & t_{i} = \max\left(a_{P_i}, T_{f, P_i} \times y_i + \left(\left(t_{i-1} + T_{P_{i-1}, P_i}\right)\times \left(1 - y_i\right) \right) \right)                          \\
                            & P_i \in  [\![ 1, n ]\!], \forall i [\![ 1, n ]\!]
    \end{align}
\end{subequations}

\subsection{Quantum Search}

The quantum computing paradigm has led to the development of several noteworthy algorithms.
Among them, Grover's search algorithm \cite{groverFastQuantumMechanical1996} is particularly relevant in the context of operations research,
as it enables the identification of solutions within an unstructured search space. The only requirement is to construct an oracle function that
returns $1$ for all states satisfying the problem's constraints, and $0$ otherwise.
The algorithm amplifies the amplitudes—and consequently the probabilities—of the desirable states, starting from a uniform superposition over
all possible configurations. Theoretically, Grover's algorithm achieves a success probability approaching $1$ after $O(\sqrt{N/M})$ \cite{boyerTightBoundsQuantum1998} iterations,
where $N$ is the size of the search space and $M$ is the number of valid solutions. As a result, the algorithm’s complexity is directly influenced
by the number of decision variables used in the problem formulation.

To perform function minimization, we adopt the quantum search-based optimization algorithm introduced by Dürr and Høyer \cite{durrQuantumAlgorithmFinding1999}.
The principle of the algorithm is as follows:
\begin{enumerate}
    \item Set an initial threshold value $k$ and add a constraint to the Grover search that only considers solutions with a cost less than $k$.
    \item If a valid solution is found, update $k$ to this solution’s cost.
    \item Repeat the process until no better solution can be found.
\end{enumerate}

We define the oracle function $O$ as follows :
\[
    O(\ket{P}, \ket{y}, k) =
    \begin{cases}
        1 & \text{if } (l(\ket{P}, \ket{y}) < k) \\ & \phantom{\text{if }} \land  (P_i \neq P_j, ~ \forall (i, j)) \\ & \phantom{\text{if }} \land (c_i \leq C^{max}, ~ \forall i \in V) \\ & \phantom{\text{if }} \land ( t_{i} \leq b_{P_i}, ~ \forall i \in V) \\
        0 & \text{otherwise}                     \\
    \end{cases}
\]

where $l$ computes the objective function \ref{eq:obj}.

Let $\hat{O}_k$ denote the oracle operator for a given threshold $k$,
and let $\hat{S}_{\ket{\Psi}}$ be the amplitude amplification operator introduced by Grover. Then, the Grover operator
is defined as $\hat{G} = \hat{O}_k \hat{S}_{\ket{\Psi}}$. Applying $\hat{G}$ iteratively amplifies the amplitudes of all
states that satisfy the constraints given in Equations~(\ref{eq:quniconst}), (\ref{eq:qcapconst}), and (\ref{eq:qtimeconst}).

Algorithm \ref{alg:global} outlines the procedure used to implement and solve the model.
We define $\ket{AD}$ as the register encoding the \textit{All Different} constraint, $\ket{\delta}$ as the qubit for the distance constraint,
and $\ket{\kappa}$ and $\ket{\tau}$ as the vectors encoding the capacity and time constraints, respectively.
For instance, $\ket{\kappa_i} = 1$ if the accumulated load of the vehicle since the start of the subtour is valid:
it equals 0 if $P_i$ is the first node of a subtour, or $\ket{c_{i-1}}$ plus the demand at the $i^{th}$ node,
provided this sum does not exceed the vehicle’s capacity.

    {\footnotesize
        \begin{algorithm}[H]
            \KwIn{$q$ the quantity vector, $T$ the travelling time matrix, $D$ the distance matrix, $k$ the cost threshold}
            \KwOut{$P$ and $y$ quantum register}
            $\ket{P} \gets H^{\otimes n \log n}$ \\
            $\ket{y} \gets H^{\otimes n}$ \\
            $\ket{q} \gets \ket{-}$ \\
            $m \gets 1$ \\
            \While {$m \leq \sqrt{2^{n\log n + n}}$}{
                \For{$i \gets 1$ \KwTo $n$}{
                    \For{$j \gets i+1$ \KwTo $n$}{
                        $\ket{AD} \gets P_i \neq P_j$
                    }
                    $\ket{c_i} \gets q_i$ \label{lalg:qinc} \\
                    \If{$\neg \ket{y_{i-1}}$}{
                        $\ket{w} \gets \ket{w} + D_{P_{i-1}, P_{i}}$ \\
                        $\ket{c_i} \gets \ket{c_{i-1}} + \ket{c_i}$ \label{lalg:c1} \\
                        $\ket{t_i} \gets \max(a_{P_i}, \ket{t_{i-1}} + T_{P_{i-1}, P_i})$ \label{lalg:t1}
                    }
                    \Else{
                        $\ket{w} \gets \ket{w} + D_{P_{i-1}, P_{depot}} + D_{P_{depot}, P_{i}}$ \\
                        $\ket{t_i} \gets \max(a_{P_i}, T_{depot, P_i})$ \label{lalg:t2}
                    }
                    $\ket{\kappa_i} \gets (\ket{c_i} \leq C^{max})$ \label{lalg:cfinal} \\
                    $\ket{\tau_i} \gets (\ket{t_i} \leq b_{P_i})$ \label{lalg:tfinal}\\
                }
                $\ket{\delta} \gets (\ket{w} \leq k)$ \\
                $\ket{q} \gets \ket{AD} \land \bigwedge_{i=0}^{n} \ket{\kappa_i} \land \bigwedge_{i=0}^{n} \ket{\tau_i} \land \ket{\delta}$ \\
                Apply Grover operator to $\ket{P} \otimes \ket{y}$ \\
                $m \gets m +1$

            }
            \caption{Grover Search CVRPTW Algorithm}
            \label{alg:global}
        \end{algorithm}
    }

\section{Implementation}

In this section, we will discuss how to implement such algorithm on a gate based quantum computer. We assume that \textit{multi-controlled X} gates also known as MCX gates are
available.

\paragraph{P register} There are various ways to represent a giant tour, which ultimately corresponds to a solution of the TSP.
The implementation of constraints depends heavily on the chosen representation. For example, \cite{srinivasanEfficientQuantumAlgorithm2018}
introduces a successor-based representation that enables efficient cost evaluation. However, as noted in \cite{zhuRealizableGASbasedQuantum2022},
it requires additional transformations to ensure that only valid solutions are marked.

In this work, we adopt a basic position-based representation, where $V \subset \mathbb{N}$ and $P_i \in V ~\forall i$. The uniqueness constraint can
be easily enforced by ensuring that each register holds a distinct value.

\subsection{Capacity and time constraints}

\paragraph{Comparisons}
As described in the previous section, both the capacity and time window constraints ultimately reduce to comparisons,
which can be efficiently implemented using a single ancilla qubit and components from the ripple-carry adder circuit introduced in \cite{cuccaroNewQuantumRipplecarry2004a}.

\paragraph{Capacity constraint}
As discussed in Section~\ref{sec:qsplit}, each value $c_i$ representing the cumulative load at node $P_i$ must be stored,
which requires $n$ dedicated quantum registers. To compute each $c_i$ as defined in Equation~\ref{eq:qcapdef},
we must first encode the demand associated with node $P_i$. For this, we use a \textit{conditional encoder},
a type of circuit that conditionally superposes multiple values onto a quantum register based on the value of one or more index registers.
An implementation of such an encoder is described in paragraph 3 of the \textsc{Methods} section in \cite{zhuRealizableGASbasedQuantum2022}, it is denoted as $F$ gate.

\begin{figure}
    \begin{center}
    \begin{quantikz}[row sep={0.6cm,between origins}, column sep=0.15cm]
        \lstick{$\ket{P}$} &\ctrl{3}& & & & & & &\\
        \lstick{$\ket{y}$} & & \gate{X} & \ctrl{2}  & & \gate{X}&& & \\
        \lstick{$\ket{c_{i-1}}$} & & & \gate[2]{Adder} & & & & &\\
        \lstick{$\ket{c_i}$} & \gate{F} & &  & & &  & \gate{\leq C^{max}} \vqw{1} && \\
        \lstick{$\ket{\kappa}$} & & & & & & & \targ{} & &

    \end{quantikz}
\end{center}
    \caption{Quantum circuit for capacity constraint}
    \label{circ:capa}
\end{figure}
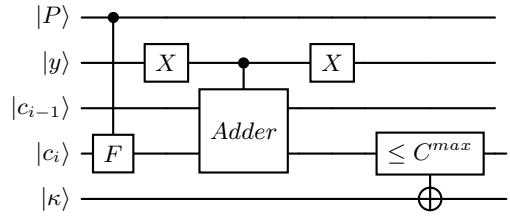

Depending on the value of the register $\ket{P_i}$, we use the encoder to assign the corresponding demand value to a load register.
If $y_i = 0$, indicating that $P_i$ is not the start of a new sub-tour, the previous accumulated load is added to this demand value.
Conversely, if $y_i = 1$, the vehicle is assumed to have returned to the depot and restarted its route, so the cumulative load is reset.
Figure \ref{circ:capa} illustrates this implementation, which directly corresponds to lines \ref{lalg:qinc}, \ref{lalg:c1}, and \ref{lalg:cfinal} of Algorithm \ref{alg:global}.

\begin{figure}
    \begin{center}
    \begin{quantikz}[row sep={0.6cm,between origins}, column sep=0.1cm]
        \lstick{$\ket{P}$}  & & \ctrl{3} & & & \ctrl{3} & & & & & & &  \\
        \lstick{$\ket{y}$} & \gate{X} & \ctrl{2} & \ctrl{1} & \gate{X} & \ctrl{2} & & &&  & & & \\
        \lstick{$\ket{t_{i-1}}$}  & & & \gate[2]{Adder} & & & & &&  & & & \\
        \lstick{$\ket{t_{i}}$} & & \gate{F}& & & \gate{F} & \gate{\leq a_{P_i }} \vqw{1} & &  \gate{Enc(a_{P_i})}  & & \gate{\leq b_{P_i}} \vqw{2} & &  \\
        \lstick{$\ket{\theta_i}$} & & & & & & \targ{} & \ctrl{1} & \ctrl{-1} &\gate{X} & \ctrl{1} & \gate{X}& \\
        \lstick{$\ket{\tau}$} & & & & & &  & \targ{} & & & \targ{} & &
    \end{quantikz}
\end{center}
    \caption{Quantum circuit for time constraint}
    \label{circ:time}
\end{figure}
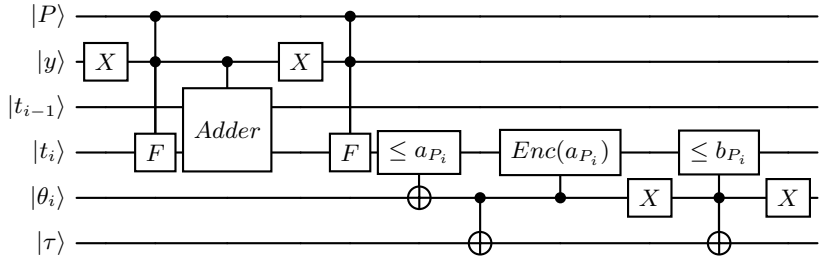

\paragraph{Time window constraint}
The time window constraint is structurally similar to the capacity constraint, with the key difference being the need to compute
a $\max()$ function. This is required to determine whether the vehicle arrives within the allowed time window.
The arrival time at node $P_i$ is computed by comparing two values: the earliest allowed delivery time $a_{P_i}$ and the current travel time.
The maximum of these two values is selected as the effective arrival time. This comparison and assignment can be implemented using a basic
comparator.
Figure \ref{circ:time} depicts this implementation, which directly maps to lines \ref{lalg:t1}, \ref{lalg:t2}, and \ref{lalg:tfinal} of Algorithm \ref{alg:global}.
The $Enc()$ gate denotes the encoding of a value into a register.

\subsection{Complexity}

\paragraph{Space complexity}
In this paper, we adopt a position-based encoding of the solution, where customer $P_i$ is interpreted as the $i^{\text{th}}$
customer visited in the giant tour. This representation requires $n \log n$ qubits to encode the tour sequence.
Additionally, we introduce $n$ qubits to store the binary variables $y_i$, which represent the splitting
decisions.

A significant number of qubits is also required to enforce the constraints. Specifically, $\mathcal{O}(n^2)$ qubits are needed
to guarantee tour validity, since there are $n(n-1)/2$ pairwise inequality constraints to ensure that each customer appears
exactly once.

For the capacity constraints, the cumulative load must be computed and stored at each step along the tour.
This requires $n (\log d_{\max} + 1)$ qubits, where $d_{\max}$ denotes the maximum vehicle load capacity.
In addition, $n$ extra qubits are needed to store boolean results indicating whether the load at each step violates
the capacity constraint. The total qubit requirement for capacity is therefore:
$n (\log d_{\max} + 1) + n = \mathcal{O}(2n + n \log d_{\max})$.

Regarding the time window constraints, a similar number of registers is required. The arrival time at each node must be
stored using $n$ registers of size $\log t_{\max} + 1$, where $t_{\max}$ is the maximum allowed time. Since the effective
arrival time is computed as the maximum between the accumulated travel time and the customer’s earliest delivery time,
$n$ additional boolean qubits are used to indicate which value is selected. Moreover, $n$ further boolean qubits are needed
to store whether the time window is respected. The total number of qubits for time-related constraints is thus:
$n (\log t_{\max} + 1) + 2n = \mathcal{O}(3n + n \log t_{\max})$.

Finally, the cost evaluation requires $\log w_{max}$ where $w_{\max}$ represents an upper bound on the solution's cost.

Putting all components together, the total number of qubits required by the algorithm is: $ \mathcal{O}(n^2 + n \log n + 6n + n \log d_{\max} + n \log t_{\max} + \log w_{max})$.
However, the number of decision variables remains relatively small, namely $n \log n + n$, as only $n$ qubits are added compared to the most compact exact quantum formulation for the TSP \cite{srinivasanEfficientQuantumAlgorithm2018}.
For instance, the example presented in \ref{example} would require at least 147 qubits to solve.

To the best of our knowledge, no oracle-based formulations of the CVRPTW exist for comparison. A few QUBO-based formulations have been proposed: route-based formulations, where the number of routes grows exponentially and must therefore
be generated using heuristics \cite{lucasIsingFormulationsMany2014}, and time-expanded formulations \cite{vargasSolvingCapacitatedVehicle2024}, where the decision variables represent every possible arc for each vehicle between nodes at different times.
In the former case, the number of decision variables depends on the number of routes considered, while in the latter it is significantly larger than in our formulations, even before accounting for potential slack variables.
In classical computing, the state of the art relies on route-based formulations,
where routes generation is handled via column generation \cite{pessoaGenericExactSolver2020}.

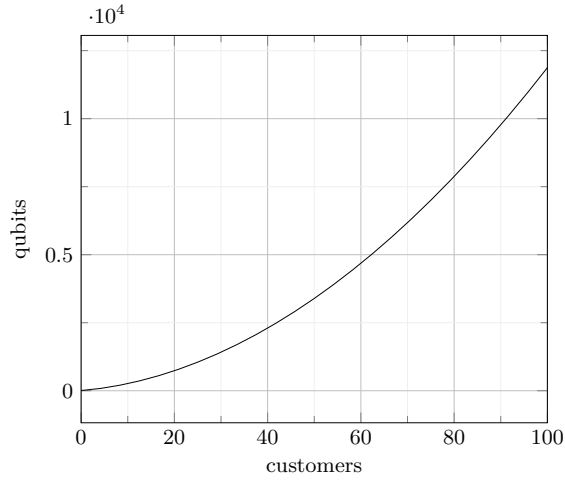
\begin{figure}[h]
    \begin{center}
        \begin{tikzpicture}[scale=0.9]
            \begin{axis}[xmin = 0,
                    xmax = 100,
                    grid = both,
                    minor tick num = 1,
                    major grid style = {lightgray},
                    minor grid style = {lightgray!25},
                    xlabel={customers},
                    ylabel={qubits}]
                \addplot[domain = 0:100] {x^2 + x*log2(x) + 6*x + x*log2(8) + x*log2(8) + log2(512)};
            \end{axis}
        \end{tikzpicture}
    \end{center}
    \caption{Number of qubits relative to the number of customers for the oracle with the maximum capacity and number of time windows both fixed to 8.}
\end{figure}

\paragraph{Time complexity}
We use the number of MCX gates as a unit for time complexity analysis.
For each gate set, we do not account for the inverse (uncomputation) circuits required to \textit{free} ancilla qubits,
as they mirror the forward operations and are already covered by the $\mathcal{O}$ notation. There are three categories of constraints contributing to the overall circuit depth.

First, the \textit{all-different} constraint, ensuring uniqueness of the tour, can be implemented using $\mathcal{O}(n^2)$ MCX gates.

Second, the \textit{capacity constraint} involves several components. The conditional encoder requires $\mathcal{O}(n \log d_{\max})$ MCX gates.
We use Cuccaro et al.'s adder \cite{cuccaroNewQuantumRipplecarry2004a}, which is efficient for this representation and contributes $\mathcal{O}(6 \log d_{\max})$ MCX gates.
The comparison circuit can also be implemented using a modified version of the same adder. Since these operations are applied at each node,
the total cost for capacity-related operations becomes: $\mathcal{O}\left(n^2 \log d_{\max} + 12n \log d_{\max}\right)$.

Third, the \textit{time window constraint} has a similar structure to the capacity constraint,
with the difference that both the encoder and the comparison are applied \textit{twice} per node.
This results in a total cost of: $\mathcal{O}\left(2n^2 \log d_{\max} + 24n \log d_{\max}\right)$.

Finally, for cost evaluation, we require $n$ \textit{conditional matrix encoders} and adders to compute the tour cost.
These operations contribute a total of: $\mathcal{O}(n^3 \log n + 6n \log n)$.
This component is relatively expensive and could potentially be optimized using a phase-based cost encoding via the Quantum Fourier
Transform (QFT). However, such an approach would require translating our representation into a form compatible with QFT-based evaluation,
which would in turn increase both qubit requirements and gate complexity \cite{zhuRealizableGASbasedQuantum2022}.
\par    In order to achieve a success probability close to 1, the Grover operator $\hat{G}$ must be applied $\mathcal{O}(\sqrt{N})$ times,
where $N$ is the size of the search space. In our case, $N = 2^{n \log n + n}$, which dominates the overall time complexity.

\section*{Conclusion}

The proposed model allows one to tackle a complex routing problem, namely the CVRPTW.
By adopting a split-inspired modeling approach on quantum computers, the model introduces only a linear additive overhead in the number of decision qubits compared
to standard quantum TSP encodings while being optimal. Consequently, the time complexity is $\mathcal{O}(\sqrt{2^{n \log n + n}})$.
Although large-scale implementations remain limited by current quantum hardware due to the gate complexity,
this research lays the groundwork for future advancements in quantum logistics optimization.


\bibliographystyle{splncs04}
\bibliography{biblio.bib}

\end{document}